\documentclass[12pt]{amsart}
\title{Linear forms and quadratic uniformity for functions on $\F_p^n$}
\author{W.T. Gowers}
\address{University of Cambridge, Department of Pure Mathematics and Mathematical Statistics, 
Wilberforce Road, Cambridge CB3 0WB, UK.}
\email{w.t.gowers@dpmms.cam.ac.uk}
\author{J. Wolf}
\address{Rutgers The State University of New Jersey, Department of Mathematics, 110 Frelinghuysen Rd., Piscataway, NJ 08854, U.S.A.}
\email{julia.wolf@cantab.net}

\usepackage{amsmath, wrapfig}
\usepackage{dsfont, a4wide, amsthm, amssymb, amsfonts, graphicx}
\usepackage{fancyhdr, xspace, psfrag, setspace, supertabular, color}

\newtheorem{theorem}{Theorem}[section]
\newtheorem{proposition}[theorem]{Proposition}
\newtheorem{lemma}[theorem]{Lemma}

\newtheorem{corollary}[theorem]{Corollary}

\newtheorem{definition}[theorem]{Definition}

\def\eps{\epsilon}
\def\E{\mathbb{E}}
\def\Z{\mathbb{Z}}
\def\R{\mathbb{R}}

\def\C{\mathbb{C}}

\def\F{\mathbb{F}}

\def\x{x}

\def\b{\beta}

\def\d{\delta}

\def\ra{\rightarrow}
\def\seq#1#2{#1_1,\dots,#1_#2}
\def\sm#1#2{\sum_{#1=1}^#2}
\def\sp#1{\langle #1\rangle}
\def\ol{\overline}
\def\hf{\hat{f}}

\begin{document}

\begin{abstract}
We give improved bounds for our theorem in \cite{Gowers:2007tcs}, which shows that a system of linear forms on $\F_p^n$ with squares that are linearly independent has the expected number of solutions in any linearly uniform subset of $\F_p^n$. While in \cite{Gowers:2007tcs} the dependence between the uniformity of the set and the resulting error in the average over the linear system was of tower type, we now obtain a doubly exponential relation between the two parameters. 

Instead of the structure theorem for bounded functions due to Green and Tao \cite{Green:2008py}, we use the Hahn-Banach theorem to decompose the function into a quadratically structured plus a quadratically uniform part. This new decomposition makes more efficient use of the $U^3$ inverse theorem \cite{Green:2008py}.

\end{abstract}

\maketitle
\tableofcontents
\section{Introduction}

In \cite{Gowers:2007tcs} we asked which systems of linear equations have the property that one can guarantee that any uniform subset of $\F_{p}^{n}$ contains the ``expected'' number of solutions. By the ``expected'' number we mean the number of solutions one would expect in a random subset of the same density, and by a ``uniform" subset of $\F_p^n$ we mean a set $A$ of density $\d$ such that the function $f_A(x)=1_A(x)-\d$ has small $U^2$ norm. 

There turns out to be a clean characterization of such systems: the number of solutions can be controlled by the $U^2$ norm in the above sense whenever the system of linear forms $\seq L m$ is \emph{square independent}, by which we mean that the functions $\seq {L^2} m$ are linearly independent over $\F_p$. The main result of \cite{Gowers:2007tcs} was the following.


\begin{theorem}\label{oldresult}
Let $L_1,\dots,L_m$ be a square-independent system of linear forms in $d$ variables of Cauchy-Schwarz complexity at most 2. For every $\epsilon>0$ there exists $c>0$ such that $f:\F_p^n\rightarrow[-1,1]$ is any function with $\|f\|_{U^2} \leq c$, then 
\[\left|\E_{\x\in(\F_p^n)^d}\prod_{i=1}^m f(L_i(\x))\right|\leq\epsilon.\]
\end{theorem}

The statement about the number of solutions in a uniform set $A \subseteq \F_p^n$ of density $\alpha$ can be recovered by setting $f$ equal to the ``balanced function" $f_A$ defined above.

We encourage the reader to consult the introduction of \cite{Gowers:2007tcs} for a detailed discussion of the context of this result, and to ignore the additional assumption of \emph{Cauchy-Schwarz complexity 2} in Theorem \ref{oldresult} for the moment. We shall not define the term here but simply remark that it is a straightforward condition that allows us to say that the average under consideration is stable under small perturbations in the $U^3$ norm. 

In \cite{Gowers:2007tcs} we defined a linear system $\seq L m$ to have \emph{true complexity} $k$ if $k$ is the least integer such that the $U^{k+1}$ norm controls the average
\[\E_{\x\in(\F_p^n)^d}\prod_{i=1}^m f(L_i(\x)).\]
It is not difficult to see that for any square-dependent system one can construct a uniform set $A \subseteq \F_p^n$ that contains significantly more than the expected number of solutions to the linear system. The proof is a generalization of the well known example of a uniform set that contains significantly more than the expected number of 4-term arithmetic progressions, which is based on the identity
\[x^2 -3(x+d)^2+3(x+2d)^2-(x+3d)^2=0.\]
Combining this fact with Theorem \ref{oldresult}, we obtain the result that a square-independent linear system has true complexity equal to 1.

We conjectured in \cite{Gowers:2007tcs} that the true complexity was always equal to the least integer $k$ such that the functions $\seq {L^k} m$ are linearly independent over $\F_p$. 

This paper is the first in a series of three in which we expand on our result in \cite{Gowers:2007tcs}. In the second paper \cite{Gowers:2009lfuII} we resolve the above conjecture, in a qualitative sense at least, for all (reasonable) systems of linear equations over $\F_p^n$, and in the third paper \cite{Gowers:2009lfuIII} we extend Theorem \ref{oldresult} to the technically more challenging setting of $\Z_N$.

In the present paper we derive a significant improvement over the bounds for Theorem \ref{oldresult} obtained in \cite{Gowers:2007tcs}, which were of tower type. Here we obtain the uniformity parameter $c$ as a doubly exponential function in the error $\eps$. 

\begin{theorem}\label{impbound}
In Theorem \ref{oldresult}, the uniformity parameter $c$ can be taken to be 
$$\exp(-\exp(c_{m,p}\eps^{-(4C_p)^{m}})),$$ where $c_{m,p}$ is a constant that depends on $m$ and $p$ only, and $C_p$ is a constant depending on $p$ only and arising in the $U^3$ inverse theorem.
\end{theorem}

The quantitative improvement which Theorem \ref{impbound} provides over Theorem \ref{oldresult} is based on a new type of decomposition of a bounded function into quadratic phases. Instead of the structure theorem for bounded functions due to Green and Tao \cite{Green:2008py}, we use the Hahn-Banach theorem to decompose the function into a quadratically structured plus a quadratically uniform part. 

This new decomposition makes more efficient use of the $U^3$ inverse theorem. Moreover, it provides a model for our more difficult proof in the cyclic group $\Z_N$ \cite{Gowers:2009lfuIII}. In this respect, we are following a course that has been strongly advocated by Green \cite{Green:2005ffm}. In that paper, we also obtain a doubly exponential bound, by following the arguments in this paper as closely as we can, but replacing subspaces with the technology of regular Bohr sets that originated in the work of Bourgain.

Just before we submitted this paper, Green and Tao \cite{Green:2010ar} proved the $\Z_N$ case of the full conjecture. Consequently the problem is now, at least in a qualitative sense, completely solved in both settings. However, the main point of this paper is the strong bounds we obtain (since we have already proved the result with a much worse bound).


\section{A simple decomposition into quadratic phases}

As in our previous paper \cite{Gowers:2007tcs}, our starting point will be the following inverse theorem of Green and Tao \cite{Green:2009ve} (in the case $p>2$) and Samorodnitsky \cite{Samorodnitsky:2007ldt} (when $p=2$).


\begin{theorem}\label{u3inverse}
Let $0<\delta\leq 1$ and let $f:\F_p^n\rightarrow\C$ be a function 
with $\|f\|_\infty\leq 1$ and $\|f\|_{U^3}\geq\delta$. Then there
exists a quadratic form $q:\F_p^n\rightarrow\F_p$ such that
$$|\E_x f(x)\omega^{q(x)}|\geq \exp(-C_p\delta^{-C_p}).$$
Here, $C_p$ is a constant that depends on $p$ only.
\end{theorem}

Green and Tao use the above theorem to decompose an arbitrary
function $f$ into two parts $f_1$ and $f_2$, where $f_2$ is
quadratically uniform and $f_1$ is quadratically structured,
in the sense that one can partition $\F_p^n$ into a small
number of quadratic subvarieties on each of which $f_1$ is
constant. In this paper, we shall take a somewhat different
approach, more closely analogous to the way conventional 
Fourier analysis is used to prove Roth's theorem. That is,
we shall simply decompose $f$ into a sum of functions of the
form $\omega^{q_i}$, where the $q_i$ are quadratic forms,
plus an error that we can afford to ignore, and then calculate 
directly using this expansion of $f$.

A big difference between the expansion we shall obtain and the
expansion of a function into Fourier coefficients is that there
does not seem to be a canonical way of doing it, because there
are far more than $p^n$ different functions of the form $\omega^q$.
(In harmonic-analysis terms, we are dealing with an ``overdetermined''
system.) This creates difficulties, which Green and Tao dealt with 
by projecting onto ``quadratic factors''. Here we
shall deal with them by applying the Hahn-Banach theorem for
finite-dimensional normed spaces.

It turns out that the Hahn-Banach theorem is a very useful tool in additive combinatorics that can be used to prove a large variety of decomposition and approximate structure theorems. It also yields a simplified proof of Green and Tao's \emph{transference principle}, which was a crucial ingredient in the proof that there exist arbitrarily long arithmetic progressions in the primes \cite{Green:2004pca}. These results are discussed at length in \cite{Gowers:2008das}.

Before we can explain why the Hahn-Banach theorem is useful, we
must state both it and one or two other simple results about
duality in normed spaces. Throughout this section we
shall refer to an inner product, which is just the standard inner 
product on $\C^n$ (or later $\C^{\F_p^n}$).


\begin{theorem}\label{hb}
Let $X=(\C^n,\|.\|)$ be a normed space and let $x\in X$
be a vector with $\|x\|\geq 1$. Then there is a vector $z$
such that $|\langle x,z\rangle|\geq 1$ and such that 
$|\langle y,z\rangle|\leq 1$ whenever $\|y\|\leq 1$.
\end{theorem}

Apart from Theorem \ref{hb}, a proof of which can be found in any standard text on functional analysis, we have aimed to keep this paper completely self-contained. We start by recalling some standard notions from the theory of normed spaces. The dual norm $\|.\|^*$ of a norm $\|.\|$ on 
$\C^n$ is defined by the formula 
\begin{equation*}
\|z\|^*=\sup\{|\langle x,z \rangle|:\|x\|\leq 1\}\\
\end{equation*}
For technical reasons, we shall generalize this concept
to the situation where the norm $\|.\|$ is defined on 
a subspace $V$ of $\C^n$. Then the dual is a seminorm,
given by the formula
\begin{equation*}
\|z\|^*=\sup\{|\langle x,z \rangle|:x\in V,\|x\|\leq 1\}\\
\end{equation*}
The next lemma is a standard fact in Banach space theory.

\begin{lemma}\label{maxdual}
Let $k$ be a positive integer, and for each $i$ between $1$ 
and $k$ let $\|.\|_i$ be a norm defined on a subspace
$V_i$ of $\C^n$. Suppose that $V_1+\dots+V_k=\C^n$, and define 
a norm $\|.\|$ on $\C^n$ by the formula
\begin{equation*}
\|x\|=\inf\{\|x_1\|_1+\dots+\|x_k\|_k:x_1+\dots+x_k=x\}
\end{equation*}
Then this formula does indeed define a norm, and its
dual norm $\|.\|^*$ is given by the formula
\begin{equation*}
\|z\|^*=\max\{\|z\|_1^*,\dots,\|z\|_k^*\}
\end{equation*}
\end{lemma}

\begin{proof}
It is a simple exercise to check that the expression does
indeed define a norm. 

Let us begin by supposing that $\|z\|_i^*\geq 1$ for some $i$.
Then there exists $x\in V_i$ such that $\|x\|_i\leq 1$ and
$|\langle x,z\rangle|\geq 1$. But then $\|x\|\leq 1$ as
well, from which it follows that $\|z\|^*\geq 1$. Therefore,
$\|z\|^*$ is at least the maximum of the $\|z\|_i^*$.

Now let us suppose that $\|z\|^*>1$.
This means that there exists $x$ such that $\|x\|\leq 1$
and $|\langle x,z\rangle|\geq 1+\epsilon$ for some
$\epsilon>0$. Let us choose $x_1,\dots,x_k$ such that
$x_i\in V_i$ for each $i$, $x_1+\dots+x_k=x$, and 
$\|x_1\|_1+\dots+\|x_k\|_k<1+\epsilon$.
Then 
\begin{equation*}
\sum_i|\langle x_i,z\rangle|>\|x_1\|_1+\dots+\|x_k\|_k
\end{equation*}
so there must exist $i$ such that 
$|\langle x_i,z\rangle|>\|x_i\|_i$, from which it follows
that $\|z\|_i^*>1$. This proves that $\|z\|^*$ is at most
the maximum of the $\|z\|_i^*$. 
\end{proof}


\begin{corollary}\label{hbcor}
Let $k$ be a positive integer and for each $i\leq k$ let 
$\|.\|_i$ be a norm defined on a subspace $V_i$ of $\C^n$,
and suppose that $V_1+\dots+V_k=\C^n$. Let 
$\alpha_1,\dots,\alpha_k$ be positive real numbers, and suppose 
that it is not possible to write the vector $x$ as a linear
sum $x_1+\dots+x_k$ in such a way that $x_i\in V_i$ for each $i$
and $\alpha_1\|x_1\|_1+\dots+\alpha_k\|x_k\|_k\leq 1$. Then 
there exists a vector $z\in\C$ such that $|\langle x,z\rangle|\geq 1$
and such that $\|z\|_i^*\leq\alpha_i$ for every $i$---or
equivalently, $|\langle y,z\rangle|\leq\alpha_i$ for every 
$i$ and every $y\in V_i$ with $\|y\|_i\leq 1$.
\end{corollary}

\begin{proof}
Let us define a norm $\|.\|$ by the formula
\begin{equation*}
\|x\|=\inf\{\alpha_1\|x_1\|_1+\dots+\alpha_k\|x_k\|_k:x_1+\dots+x_k=x\}
\end{equation*}
Then our hypothesis is that $\|x\|\geq 1$. Therefore, by Theorem
\ref{hb} there is a vector $z$ such that $|\langle x,z\rangle|\geq 1$
and $|\langle y,z\rangle|\leq 1$ whenever $\|y\|\leq 1$. 

The second condition tells us that $\|z\|^*\leq 1$, and Lemma
\ref{maxdual}, applied to the norms $\alpha_i\|.\|_i$, tells us 
that $\|z\|^*$ is the maximum of the numbers $\alpha_i^{-1}\|z\|_i^*$.
Therefore, $\|z\|_i^*\leq\alpha_i$ for every $i$, as stated.
\end{proof}

Recall that the difficulty we are trying to deal with is that there is
no (known) canonical way of decomposing a function into functions of
the form $\omega^q$. Corollary \ref{hbcor} is extremely helpful 
for proving the existence of decompositions under these circumstances.
Instead of trying to find a decomposition explicitly, one assumes that
there is no such decomposition and uses Corollary \ref{hbcor} to derive a
contradiction. The next result illustrates the technique.


\begin{theorem}\label{quadfourier}
Let $f:\F_p^n\rightarrow\C$ be a function such that 
$\|f\|_2\leq 1$. Then for every $\delta>0$ and $\eta>0$
there exists $M$ such that $f$ has a decomposition of
the form
\begin{equation*}
f(x)=\sum_i\lambda_i\omega^{q_i(x)}+g(x)+h(x),
\end{equation*}
where the $q_i$ are quadratic forms on $\F_p^n$, and
\begin{equation*}
\eta^{-1}\|g\|_1+\delta^{-1}\|h\|_{U^3}+M^{-1}\sum_i|\lambda_i|\leq 1.
\end{equation*}
In fact, $M$ can be taken to be $\exp(C_p(\eta\delta)^{-C_p})$, where $C_p$ is the constant in Theorem \ref{u3inverse}.
\end{theorem}

\begin{proof}
Suppose not. Then for every quadratic form $q$ on $\F_p^n$
let $V(q)$ be the one-dimensional subspace of $\C^{\F_p^n}$
generated by the function $\omega^q$, with the obvious
norm: the norm of $\lambda\omega^q$ is~$|\lambda|$.

Applying Corollary \ref{hbcor} to these norms and subspaces
and also to the $L_1$  norm and $U^3$ norm defined on all 
of $\C^{\F_p^n}$, we deduce that there is a function
$\phi:\F_p^n\rightarrow\C$ such that $|\langle f,\phi\rangle|\geq 1$,
$\|\phi\|_\infty\leq\eta^{-1}$, $\|\phi\|_{U^3}^*\leq\delta^{-1}$
and $|\langle \phi,\omega^q\rangle|\leq M^{-1}$ for every 
quadratic form $q$. 

Now the fact that $|\langle f,\phi\rangle|\geq 1$ implies, by
Cauchy-Schwarz, that $\|\phi\|_2\geq 1$. But we also know
that $\langle \phi,\phi\rangle\leq \|\phi\|_{U^3}\|\phi\|_{U^3}^*$,
so $\|\phi\|_{U^3}\geq\delta$. Applying the inverse theorem to
$\eta\phi$, we find that there is a quadratic form $q$ such
that $|\langle \phi,\omega^q\rangle|\geq \exp(-C_p(\eta\delta)^{-C_p})$,
contradicting the fact that it has to be at most $M^{-1}$.
\end{proof}

Just before we continue, let us briefly discuss a more obvious
approach to a slight variant of Theorem \ref{quadfourier} and see why it
does not work.  Theorem \ref{u3inverse} tells us that every bounded
function $f$ with large $U^3$ norm correlates well with some function
of the form $\omega^q$. So one might assume that $f$ is bounded and
try a simple inductive argument
along the following lines. If $\|f\|_{U^3}$ is large, then Theorem
\ref{quadfourier} gives us a quadratic form $q_1$ such that $f$
correlates with $\omega^{q_1}$. So choose $\lambda_1$ such that
$\|f-\lambda_1\omega^{q_1}\|_2$ is minimized, and let
$f_1=f-\lambda_1\omega^{q_1}$. Because of the correlation, $\|f_1\|_2$
is substantially less than $\|f\|_2$. Now repeat for $f_1$, and keep
going until you reach some $k$ for which $\|f_{k+1}\|_{U^3}$ is small.

The problem with this argument is that we gradually lose control of
the boundedness of $f$. As we continually subtract the functions
$\lambda_i\omega^{q_i}$ the $L_2$ norm goes down, but the $L_\infty$   norm can go up. And $L_2$ control is not enough for Theorem
\ref{u3inverse} as the example of a suitably normalized arithmetic progression shows. It turns out that a variant of the inductive argument outlined above can be made to work if one uses a weaker assumption than boundedness \cite{Candela:2007pu3} (which means that the result proved is stronger). Green and Tao's approach to quadratic Fourier analysis assumes an $L_\infty$   bound for $f$ and uses averaging
projections, which decrease both the $L_2$  and $L_\infty$ norms. Thus,
there seems to be a genuine difference between Theorem
\ref{quadfourier} and their approach.

However, there are two aspects of Theorem \ref{quadfourier} that place
considerable limits on how useful it is. The first is that $M$ is
rather large, so that bounds that depend on the theorem tend to be
rather large as well. The second, which is more serious, is that there
is no useful bound on the number of quadratic phase functions used to
decompose $f$. We shall deal with these two problems in turn.

\section{Introducing quadratic averages}

In order to reduce $M$, we shall use a slightly stronger form of
Theorem \ref{u3inverse}, which Green and Tao \cite{Green:2008py} mention but do not
need, and therefore do not formally state. To begin with, here is
a variant that they do state. If $V$ is a subspace of 
$\F_p^n$ and $y\in\F_p^n$, then they define a seminorm $\|.\|_{u^3(y+V)}$
on functions from $\F_p^n$ to $\C$ by setting 
$$\|f\|_{u^3(y+V)}=\sup_q |\E_{x\in y+V}f(x)\omega^{-q(x)}|,$$ 
where the supremum is taken over all quadratic forms $q$ on $y+V$. 


\begin{theorem} \label{betterinverse}
Let $f:\F_p^n\ra\C$ be a function such that
$\|f\|_\infty\leq 1$ and $\|f\|_{U^3}\geq\delta$. 
Then there exists a subspace $V$ of $\F_p^n$ of codimension at
most $(2/\d)^{C_p}$, where $C_p$ is a constant that depends only
on $p$, with the property that 
\[\E_y\|f\|_{u^3(y+V)}\geq(\d/2)^{C_p}.\]
\end{theorem}

One can deduce Theorem \ref{u3inverse} very
simply from this version: by an averaging argument, there must exist
$y$ such that $f$ correlates well on $y+V$ with some quadratic
phase function $\omega^q$; this function can be extended to a 
function on the whole of $\F_p^n$ in many different ways, and
a further averaging argument yields Theorem \ref{u3inverse}. 

It turns out that, as Green and Tao remark, a slightly more precise
theorem holds. The result as stated tells us that for each $y$ we
can find a local quadratic phase function $\omega^{q_y}$ defined on $y+V$
such that the average of $|\E_{x\in y+V}f(x)\omega^{q_y(x)}|$ is at
least $(\d/2)^C$. However, it is actually possible to do this in 
such a way that the ``quadratic parts'' of the quadratic phase 
functions $q_y$ are the same. More precisely, it can be done in 
such a way that each $q_y(x)$ has the form $q(x-y)+\phi_y(x-y)$ for
some quadratic function $q:V\ra\F_p$ (that is independent of $y$) and 
some linear functionals $\phi_y:V\ra\F_p$. A similar statement can be read out of \cite{Samorodnitsky:2007ldt} for the case of $\F_2^n$ (see also \cite{Wolf:2009lit}).

This will be convenient to us later, so let us make a definition
so that we can refer to this property concisely.
 
\begin{definition}\label{defquadave} Let $V$ be a subspace of $\F_p^n$ and let $q$ be
a quadratic form on $V$. A \emph{quadratic average with
base} $(V,q)$ is a function of the form $Q(x)=\E_{y\in
x-V}\omega^{q_y(x)}$, where each function $q_y$ is a 
quadratic map from $y+V$ to $\F_p$ defined by a formula
of the form $q_y(x)=q(x-y)+\phi_y(x-y)$ for some Freiman homomorphism
$\phi_y:V\ra\F_p$. The \emph{rank} of $Q$ is the rank
of the quadratic form $q$, and the \emph{complexity} of $Q$
is the codimension of $V$.
\end{definition}

Notice that if $Q$ is a quadratic average then $\|Q\|_\infty\leq 1$. 
For fixed $x$ it is natural to write the set where 
$V(x-y)$ is defined as $x-V$, but since $V$ is symmetric this is 
equal to $x+V$. We write it as $x-V$ because writing at as
$x+V$ makes certain proofs slightly confusing.

Before we move on let us describe a particular property of quadratic averages which will be used on a number of occasions in the sequel. It is a consequence of the fact that the quadratic phases are somewhat ``parallel". In order to become familiar with quadratic averages, let us prove that high-rank quadratic averages have small $U^2$ norm. We begin with a standard fact about Gauss sums. Since the proof is short, we include it for the sake of completeness.


\begin{lemma} \label{gauss}
Let $q$ be a quadratic form of rank $r$ and let $\phi$ be a linear
function. Then 
\[|\E_x\omega^{q(x)+\phi(x)}|\leq p^{-r/2}.\]
\end{lemma}

\begin{proof}
Let $\beta$ be the bilinear form given by the formula
$\beta(x,y)=q(x+y)-q(x)-q(y)$. Then
\begin{equation*}
|\E_x\omega^{q(x)+\phi(x)}|^2=\E_{x,y}\omega^{q(x+y)-q(x)+\phi(x+y)-\phi(x)}
=\E_{x,y}\omega^{q(y)+\phi(y)+\beta(x,y)}.
\end{equation*}
For each $y$, the expectation over $x$ is zero unless 
$\beta(x,y)$ is constant, and therefore identically zero,
in $x$. But for such a $y$ we have $q(y)=0$ as well, so
the expectation over $x$ is $\omega^{\phi(y)}$, which has
modulus 1.

Since $q$ has rank $r$, the space of $y$ such that $\beta(x,y)$
is zero for every $x$ has codimension $r$, so the right hand
side is at most $p^{-r}$. This proves the lemma.
\end{proof}

Note that in fact we have the more precise result that the expectation
is zero if $\phi$ does not vanish on the annihilator of $\beta$, and
$p^{-r/2}$ otherwise.


\begin{lemma}\label{quadaveu2}
Let $Q$ be a quadratic average of rank $r$. Then $\|Q\|_{U^2}\leq p^{-r/4}$.
\end{lemma}

\begin{proof}
Let $V$ be a subspace of $\F_p^n$ and let $\phi_1$, $\phi_2$, 
$\phi_3$ and $\phi_4$ be linear functions defined on $V$. 
Let $q$ be a quadratic form on $V$ of rank $r$ and let
$q_i=q+\phi_i$. Then
\[\E_{x,a,b\in V}\omega^{q_1(x)-q_2(x+a)-q_3(x+b)+q_4(x+a+b)}
=\E_{a,b\in V}\omega^{\b(a,b)},\]
where $\b$ is a bilinear form of rank $r$. For each $a$, the
expectation over $b$ is zero unless $\b(a,b)$ is zero for
every $b$, which happens only when $a$ belongs to an $r$-codimensional
subspace of $V$. Therefore, the right-hand side equals $p^{-r}$.

Now let $Q$ be a quadratic average with base $(V,q)$. Then 
\[\|Q\|_{U^2}^4=\E_{x_1+x_2=x_3+x_4}Q(x_1)Q(x_2)\ol{Q(x_3)Q(x_4)}.\]
If we condition on the right-hand side according to which translates
$x_1$, $x_2$, $x_3$ and $x_4$ belong to, we obtain an expectation
of expectations of the form discussed in the previous paragraph.
Indeed, let $y_1+y_2=y_3+y_4$ and let us calculate the average
over all $x_1+x_2=x_3+x_4$ such that $x_i\in y_i+V$. Inside $y_i+V$,
we have $Q(x)=\omega^{q(x-y_i)+\phi_{y_i}(x-y_i)}$ for some Freiman 
homomorphism $\phi_i$. Thus, the average in question is
\begin{equation*}
\E_{x_1+x_2=x_3+x_4}\omega^{q(x_1-y_1)+\phi_{y_1}(x_1-y_1)+
q(x_2-y_2)+\phi_{y_2}(x_2-y_2)-q(x_3-y_3)-\phi_{y_3}(x_3-y_3)-
q(x_4-y_4)-\phi_{y_4}(x_4-y_4)},
\end{equation*}
where the average is over $x_i\in y_i+V$. we can then substitute
by setting a new $x_i$ to be $x_i-y_i$ and we have an average of
the required form. The result follows.
\end{proof}

It may seem slightly strange that the above result does not depend
on the complexity of the quadratic averages $Q$. But
this is because the assumption about the rank of $q$ is stronger
the larger the codimension of $V$. Thus, there is in fact a 
dependence but it is disguised by the way the lemma is formulated.

Now we can state a slightly more precise version of the inverse
theorem that follows easily from the remark of Green and Tao when $p>2$, and 
follows with some care from the work of Samorodnitsky in the case $p=2$ \cite{Wolf:2009lit}.


\begin{theorem}\label{fancyinverse}
Let $f:\F_p^n\ra\C$ be a function such that $\|f\|_\infty\leq 1$ and
$\|f\|_{U^3}\geq\delta$. Then there exists a quadratic average $Q$
of complexity at most $(2/\d)^{C_p}$ such that \[|\sp{f,Q}|\geq(\d/2)^{C_p}/2.\]
\end{theorem}

\begin{proof}
The results of Green and Tao tell us that we can find a subspace 
$V$ satisfying the above codimension bound,
a quadratic function $q$ defined on $V$, and for each $y$ a
linear map $\phi_y:V\ra\F_p$, such that, defining
$q_y(x)=q(x-y)+\phi_y(x-y)$ on $y+V$, we have
\[\E_y|\E_{x\in y+V}f(x)\omega^{-q_y(x)}|\geq(\d/2)^{C_p}.\]
For each function $q_y$ we can add a constant $\lambda_y \in \F_p$ without
affecting the left-hand side. We can choose this constant so that
\[\Re(\E_{x\in y+V}f(x)\omega^{-q_y(x)+\lambda_y})
\geq \frac{1}{2}|\E_{x\in y+V}f(x)\omega^{-q_y(x)}|.\]
Therefore, after suitably redefining the functions $q_y$ and
setting $Q(x)=\E_{y\in x-V}\omega^{q_y(x)}$, we have
\[
|\sp{f,Q}|\geq \Re(\E_x\E_{y\in x-V}f(x)\omega^{-q_y(x)})
\geq \frac{1}{2}\E_y|\E_{x\in y+V}f(x)\omega^{-q_y(x)}|,
\]
which proves the theorem.
\end{proof}

Now let us quickly deduce a corresponding decomposition result
by the same method as before.


\begin{theorem}\label{fancyquadfourier}
Let $f:\F_p^n\rightarrow\C$ be a function such that 
$\|f\|_2\leq 1$. Then for every $\delta>0$ and $\eta>0$
there exists $M$ such that $f$ has a decomposition of
the form
\begin{equation*}
f(x)=\sum_i\lambda_iQ_i(x)+g(x)+h(x),
\end{equation*}
where the $Q_i$ are quadratic averages on $\F_p^n$ of complexity at most $(2/\delta)^{C_p}$, and
\begin{equation*}
\eta^{-1}\|g\|_1+\delta^{-1}\|h\|_{U^3}+M^{-1}\sum_i|\lambda_i|\leq 1.
\end{equation*}
In fact, $M$ can be taken to be $(2/\eta\d)^{C_p}/2$.
\end{theorem}

\begin{proof}
Suppose not. Then for every quadratic average $Q$ on $\F_p^n$,
let $V(Q)$ be the one-dimensional subspace of $\C^{\F_p^n}$
generated by $Q$ with the norm of $\lambda Q$ defined to be
$|\lambda|$.

Applying Corollary \ref{hbcor} to these norms and subspaces,
and also to the $L_1$  norm and $U^3$ norm defined on all 
of $\C^{\F_p^n}$, we find a function
$\phi:\F_p^n\rightarrow\C$ such that $|\langle f,\phi\rangle|\geq 1$,
$\|\phi\|_\infty\leq\eta^{-1}$, $\|\phi\|_{U^3}^*\leq\delta^{-1}$
and $|\langle \phi,Q\rangle|\leq M^{-1}$ for every 
quadratic average $Q$. 

As before, our assumptions imply that $\|\phi\|_{U^3}\geq\delta$. 
Applying the more precise inverse result (Theorem \ref{fancyinverse}) to
$\eta\phi$, we find that there is a quadratic average $Q$ of complexity at most $(2/\delta)^{C_p}$ such that
$|\langle\phi,Q\rangle|\geq(\eta\d/2)^{C_p}$,
contradicting the fact that this correlation has to be at most $M^{-1}$.
\end{proof}

\section{Clusters of highly correlated quadratic phases}

Now let us deal with the difficulty that the above decomposition may
be into a huge number of quadratic averages. Corollary \ref{unitvectors}
below is a simple result that will be used to place a strong restriction on
the decompositions that can occur.

To prepare for it, let us define a \textit{vertex-weighted}
graph to be a graph $G$ together with a function $\mu:V(G)\ra\R_+$. The
weights $\mu$ induce an obvious weighting of the edges: if $xy$
is an edge of $G$ then its weight is $\mu(x)\mu(y)$. Finally,
if $H$ is a subgraph of $G$, then we define the \textit{total
weight} of $H$ to be the sum of the weights of all the edges 
in $H$.


\begin{lemma} \label{nbds}
Let $C$ and $\gamma>0$ be constants, and let $G$ be a weighted graph 
with total weight at most $C$. Then there exist vertices
$x_1,\dots,x_k$ of $G$, with $k\leq C^2/2\gamma$, such that
if $A$ is the set of all vertices not joined to any of the
$x_i$, then the total weight of the subgraph of $G$ induced
by $A$ is at most $\gamma$.
\end{lemma}

\begin{proof}
If the total weight of $G$ is at least $\gamma$, then 
$\sum_x\sum_{y\in N(x)}\mu(x)\mu(y)\geq\gamma/2$, where
we write $N(x)$ for the neighbourhood of $x$. Since
$\sum_x\mu(x)\leq C$, it follows that there exists $x$
such that $\sum_{y\in N(x)}\mu(y)\geq\gamma/2C$. Let $x_1$ be
a vertex with this property and let $G_1$ be the subgraph
of $G$ induced by the vertices not in $N(x_1)$. 

Now repeat this argument for $G_1$, and so on. Since at 
each stage we remove a set of vertices with weights summing 
to at least $\gamma/2C$, the process cannot continue for more
than $C^2/2\gamma$ steps, at which point the graph induced by
the complements of all the neighbourhoods we have chosen has
total weight at most $\gamma$, as claimed.
\end{proof}


\begin{corollary} \label{unitvectors}
Let $\seq u n$ be a collection of vectors of norm at most 1 in a
Hilbert space $H$, let $\seq \lambda n$ be scalars with $\sm i
n|\lambda_i|\leq C$ and let $\d>0$. Then there are vectors
$u_{i_1},\dots,u_{i_k}$ and a set $A\subseteq\{1,2,\dots,n\}$ such that
$k\leq C^2/\d^2$, and with the following properties. For every
$i\notin A$ there exists $j$ such that
$|\sp{u_i,u_{i_j}}|\geq\d^2/2C^2$, and $\|\sum_{i\in
A}\lambda_iu_i\|_2\leq\d$.
\end{corollary}

\begin{proof}
Define a vertex-weighted graph by taking the vectors $\seq u n$
as the vertices, letting the weight of $u_i$ be $|\lambda_i|$,
and joining $u_i$ to $u_j$ if and only if
$|\sp{u_i,u_j}|\geq\d^2/2C^2$. Now we apply the previous lemma
with $\gamma=\d^2/2$. It gives us vectors $u_{i_1},\dots,u_{i_k}$, 
with $k\leq C^2/\d^2$,
such that if $A$ is the set of all $i$ for which there is
no $i_j$ with $|\sp{u_i,u_{i_j}}|\geq\d^2/2C^2$, then the sum
of all $|\lambda_i||\lambda_j|$ with $i,j\in A$ and
$|\sp{u_i,u_{i_j}}|\geq\d^2/2C^2$ is at most $\d^2/2$. But
\[\|\sum_{i\in A}\lambda_iu_i\|_2^2=\sum_{i,j\in A}|\lambda_i|
|\lambda_j||\sp{u_i,u_j}|\]
and we can split the last sum into two parts, according to 
whether $|\sp{u_i,u_j}|$ is at least $\d^2/2C^2$ or less than
$\d^2/2C^2$. The first part is at most $\d^2/2$, since we always
have $|\sp{u_i,u_j}|\leq 1$, and the second is at most 
$C^2\d^2/2C^2=\d^2/2$. This proves the result. 
\end{proof} 

We are interested in sums of the form $\sum_i\lambda_iu_i$ where the
$u_i$ are built out of polynomial phase functions. These have the
property that if two of them have a significant correlation then they
must have a strong algebraic relationship. For example, if two linear
phase functions have any correlation at all, then they must be equal
up to a scalar multiple, and if two quadratic phase functions are well
correlated, then the difference of the corresponding quadratic forms
must have low rank. Thus, Corollary \ref{unitvectors} immediately
implies that if $f_1=\sum_i\lambda_i\omega^{q_i}$, then we can write
$f_1$ as $f_2+g$, where $f_2$ is composed of a few clusters of
quadratic phase functions that do not differ except in a ``linear''
way, and $g$ is small in $L_2$.

Now let us make this remark precise. We start off by finding out what it means for two quadratic averages 
to be well correlated.


\begin{lemma} \label{quadavecorr}
Let $Q$ and $Q'$ be quadratic averages with bases $(V,q)$ and
$(V',q')$, respectively. Suppose that the rank of $q-q'$ considered
as a quadratic form on $V\cap V'$ is $r$. Then $|\sp{Q,Q'}|\leq
p^{-r/2}$.
\end{lemma}

\begin{proof}
By the remarks following Definition \ref{defquadave}, we know that on each translate of $V\cap V'$ the function $Q\ol{Q'}$ is given by a 
formula of the form $\omega^{q-q'-\phi}$, where $\phi$ is a linear
function. Therefore, by Lemma \ref{gauss}, the expectation of
$Q\ol{Q'}$ over each such translate has modulus at most $p^{-r/2}$. 
The result follows.
\end{proof}

Once again, the apparent lack of any dependence on the complexities of $Q$ and $Q'$ is an illusion: the higher the complexity, the stronger
the rank assumption.

We now examine the $U^2$ dual norm of low-rank quadratic averages. 
First we need a rather crude lemma. Given a subspace $W$ of $\F_p^n$,
we define $W^\perp$ to be the space of all $r$ such that
$r^Tx=0$ for every $x\in W$. 


\begin{lemma} \label{calculation}
Let $W$ be a subspace of $\F_p^n$ of codimension $d$, let $y\in W^\perp$
and let $\phi$ be a linear function from $W$ to $\F_p$. Let
$g(x)=\omega^{\phi(x-y)}$ when $x\in y+W$ and 0 otherwise. Then
$\|g\|_{U^2}^*=p^{-d/4}$.
\end{lemma}

\begin{proof} 
Since the $U^2$ dual norm is unaffected by translation and by
multiplying by a linear phase function, we can assume that $y=0$
and $\phi$ is the zero function. Thus, we are calculating the
$U^2$ dual norm of the characteristic function of $W$. The
Fourier transform of this function takes the value $p^{-d}$
at every $y\in W^\perp$ and 0 everywhere else. The $U^2$ dual
norm is the $\ell^{4/3}$ norm of this Fourier transform,
which is $(p^dp^{-4d/3})^{3/4}$, which equals $p^{-d/4}$.
\end{proof}

Before we continue, here is an alternative proof of
Lemma \ref{calculation} that does not use the Fourier transform.
This will be useful later when we want to generalize it. One
observes first that $g(x)=p^{2d}\E_{z+w-y=x}g(z)g(w)\ol{g(y)}$.
From this it follows that, for any function $f:\F_p\ra\C$,
\[\sp{g,f}=p^{2d}\E_{x+y=z+w}f(x)g(y)\ol{g(z)g(w)},\]
which has modulus at most $p^{2d}\|f\|_{U^2}\|g\|_{U^2}^3$.
It is easy to prove (without Fourier analysis) that
$\|g\|_{U^2}^4=p^{-3d}$, so $|\sp{f,g}|\leq p^{-d/4}\|f\|_{U^2}$,
which proves the lemma. We remark that the inequality is sharp because
$\sp{g,g}=p^{-d/4}\|g\|_{U^2}$.


\begin{lemma}\label{quadaveu2*}
Let $V$ and $V'$ be subspaces of $\F_p^n$, and let $q$ and $q'$ be
quadratic forms defined on $V$ and $V'$, respectively. Suppose that the codimension
of $V\cap V'$ is $d$ and the rank of the
restriction of $q-q'$ to $V\cap V'$ is $r$. Let $Q$ and $Q'$ be
quadratic averages with bases $(V,q)$ and $(V',q')$,
respectively. Then $\|Q\ol{Q'}\|_{U^2}^*\leq p^{3(d+r)/4}$.
\end{lemma}

\begin{proof}
For any fixed $y$, the restriction of the function $Q\ol{Q'}$ to $y+V\cap V'$ is equal to
$\omega^{q-q'+\phi}$ for some linear function $\phi$. Since the rank
of $q-q'$ is $r$, there is a subspace $W\subset V\cap V'$, of
codimension $r$ in $V\cap V'$, such that the function $\omega^{q-q'}$
is constant on those translates of $W$ that live inside $y+V\cap
V'$. Thus the restriction of $Q\ol{Q'}$ to any translate of $W$ that lives inside $y+V\cap V'$ is 
a linear phase function. By Lemma \ref{calculation}, each of these 
restrictions has $U^2$ dual norm $p^{-(r+d)/4}$, so their sum, which is
the restriction of $Q\ol{Q'}$ to $y+V\cap V'$, has $U^2$ dual norm at most 
$p^{r-(r+d)/4}$. The same is true for all translates of $V\cap V'$, of which there
are $p^d$. Therefore, the $U^2$ dual norm of $f$ is at most
$p^{3(r+d)/4}$, as claimed.
\end{proof}

Now let us put these facts together in order to obtain a more sophisticated decomposition into linear combinations of
quadratic averages.


\begin{theorem} \label{quadfouriersmallk}
Let $f:\F_p^n\ra\C$ be a function such that $\|f\|_2\leq 1$, and
let $\delta>0$. Let $d=(2/\d)^{C_p}$ and 
$C=(2/\d^2)^{C_p}$. Then $f$ has a decomposition
\begin{equation*}
f(x)=\sm i k Q_i(x)U_i(x)+g(x)+h(x),
\end{equation*}
where $k\leq C^2/\d^2$, the $Q_i$ are quadratic averages on 
$\F_p^n$, $\sm i k\|U_i\|_{U^2}^*\leq 2^{3/2}C^{4}\d^{-3}p^{3d/4}$, 
$\sm i k\|U_i\|_\infty\leq C$, $\|g\|_1\leq 2\d$ and
$\|h\|_{U^3}\leq\d$.
\end{theorem}

\begin{proof}
We begin by using Theorem \ref{fancyquadfourier} to decompose $f$ 
into a sum $\sum_i\lambda_iQ_i(x)+g'(x)+h(x)$, where each $Q_i$ is a
quadratic average of complexity at most $(2/\d)^{C_p}$, 
and $\|g'\|_1\leq\d$, $\|h\|_{U^3}\leq\d$ and 
$\sum_i|\lambda_i|\leq (2/\d^2)^{C_p}$. Let us set $d$ to equal
$(2/\d)^{C_p}$ and $C$ to equal $(2/\d^2)^{C_p}$.

Now we apply Corollary \ref{unitvectors} to the linear combination
$\sum_i\lambda_iQ_i$. Without loss of generality, the functions that
it gives us are $\seq Q k$. Hence we can write $\sum_i\lambda_iQ_i$ in the form
$\sm i kQ_iU_i+g''$, where $k\leq C^2/\d^2$, $\|g''\|_2\leq\d$ and 
each $U_i$ is a function
of the form $\sum_{j\in A_i}\lambda_j\ol{Q_i}Q_j$ with 
$|\sp{Q_i,Q_j}|\geq\d^2/2C^2$.

By Lemma \ref{quadavecorr} it follows that for each $j\in A_i$ the
rank $r$ of the quadratic form $q_i-q_j$ is such that
$p^{r/2}\leq 2C^2/\d^2$. Since the complexity of each $Q_i$ is at
most $d$, Lemma \ref{quadaveu2*} then tells us that
$\|\ol{Q_i}Q_j\|_{U^2}^*$ is at most $(2C^2/\d^2)^{3/2}p^{3d/4}$.
Since $\sm i k\sum_{j\in A_i}|\lambda_j|\leq C$, it follows that
$\sm i k\|U_i\|_{U^2}^*\leq 2^{3/2}C^{4}\d^{-3}p^{3d/4}$ and
$\sm i k\|U_i\|_\infty\leq C$. Since $\|g''\|_1\leq\|g''\|_2$, 
the result is true with $g=g'+g''$.
\end{proof}

By making some very minor changes to the argument it is possible to gain independent control over the $L_1$ and the $U^3$ error in this approximation, but we shall not need to do so here.

\section{A stronger decomposition for highly uniform functions}

Later on, we shall need to use the fact that if a function $f$
is uniform, then the quadratic averages used to decompose it
can all be taken to have high rank. This result is very 
plausible, as low-rank quadratic averages are anti-uniform
and should therefore not be necessary, but there does not
seem to be a truly short proof of this fact.

The next lemma shows that one can split any function of the type
$\sm i kQ_iU_i$ into a low-rank part and a high-rank part in 
such a way that there is a substantial gap between the ranks 
in the two parts. The proof is very short, mainly because 
the work has already been done: in order to have a useful result
we are relying on the fact that $k$ is small.


\begin{lemma}\label{rankgap}
Let $R_0$, $m$ and $t >1$ be 
constants. Let $Q_1, Q_2, \dots, Q_k$ be quadratic averages. Then there is a partition of
$\{1,2,\dots,k\}$ into two sets $L$ and $H$, and a constant
$R\in[R_0,m^k(R_0+t)]$, such that the rank of $Q_i$ is at most $R$ for every 
$i\in L$ and at least $mR+t$ for every $i\in H$.
\end{lemma}

\begin{proof}
Without loss of generality the $Q_i$ are arranged in increasing
order of rank. If there is no $i$ such that $Q_i$ has rank at
least $m^i(R_0+t)$, then let $L=\{1,\dots,k\}$, and we are done.
Otherwise, for each $j$ let $R_j=m^{j}R_{0}+(m^{j-1}+m^{j-2}+\dots+1)t$ and let $i$ be minimal such that $Q_i$ has rank at least $R_i$. Since $R_i=mR_{i-1}+t$ and $R_k\leq m^k(R_0+t)$, we are done.
\end{proof}

The next lemma is a standard application of Bogolyubov's method, and shows that anti-uniform functions are well approximated by their convolution with a low-codimensional subspace $V$. Alternatively, in the language of Green and Tao, which we shall not be using here, anti-uniform functions are well-approximated by their projection onto a suitable low-complexity linear factor. We write $\mu_V$ for the characteristic measure of this subspace $V$, which has the property that $\E_x \mu_V(x)=1$.


\begin{lemma}\label{bogolyubov}
Let $\d>0$ and $T$ be constants, let $f:\F_p^n\ra\C$ and suppose that
$\|f\|_{U^2}^*\leq T$. Then there is a linear subspace $V$ of
codimension at most $\d^{-4}T^4$ such that $\|f-f*\mu_V\|_2\leq\d$.
\end{lemma}

\begin{proof} 
Our assumption can be rephrased in Fourier terms as the assertion
that $\|\hf\|_{4/3}\leq T$. By the Fourier inversion formula, we
know that $f(x)=\sum_r \hat{f}(r)\omega^{-r^Tx}$ for every $x$. Let 
$\rho=\d^{3}T^{-2}$ and let $K=\{r:|\hf(r)|\geq\rho\}$.

Let $V$ be the subspace of all $x$ such that $r^Tx=0$ for every $r\in
K$. Since $\|\hat{f}\|_{4/3}\leq T$, we know that $|K|\leq
T^{4/3}\rho^{-4/3}$. Therefore, $V$ has codimension at most
$(T/\rho)^{4/3} =\delta^{-4}T^4$.  Then we can 
decompose $f$ as a sum $f_1+f_2$, where 
$f_1(x)=\sum_{r\in K}\hat{f}(r)\omega^{-r^Tx}$ and
$f_2(x)=\sum_{r\notin K}\hat{f}(r)\omega^{-r^Tx}$. It is easy to see that $f_1=f*\mu_V$, and we can bound the $L_2$ norm of $f_2$ as follows.
\begin{equation*}
\|f_2\|_2^2=\|\hat{f_2}\|_2^2\leq\|\hat{f_2}\|_{4/3}^{4/3}
\|\hat{f_2}\|_\infty^{2/3}\leq\rho^{2/3}T^{4/3}=\delta^{2}.
\end{equation*}
Therefore, the statement of the lemma follows from our calculations.
\end{proof}

It follows from Lemma \ref{bogolyubov} that the product of a low-rank quadratic average with an anti-uniform function is also well approximated by its convolution with a suitable subspace of bounded codimension.


\begin{corollary}\label{quadavetimesu2*}
Let $\d>0$ and $T$ be constants, let $U:\F_p^n\ra\C$ and suppose
that $\|U\|_{U^2}^*\leq T$. Let $r$ and $d$ be constants, and let 
$Q$ be a quadratic average of rank at most $r$ and complexity at
most $d$. Let $f=QU$. Then there is a linear subspace $V$ of 
codimension at most $\d^{-4}T^4+r+p^d$ such that $\|f-f*\mu_V\|_2\leq\d$. 
\end{corollary}

\begin{proof}
By Lemma \ref{bogolyubov} there is a subspace $V_0$ of
codimension at most $\d^{-4}T^4$ such that $\|U-U*\mu_{V_0}\|_2\leq\d$.
Let $Q$ have base $(V_1,q)$. Since $q$ has rank at most $r$, 
$V_1$ contains a subspace $W_1$ of codimension at most $r$ such 
that $\omega^q$ is constant on each translate of $W_1$. Therefore,
for each $y$ the function $\omega^{q(x-y)}$ is constant on 
those translates of $W_1$ that are subsets of $y+V_1$. 

Now on $y+V_1$ the quadratic average $Q(x)$ takes the form 
$\omega^{q(x-y)+\phi_y(x-y)}$ for some linear form $\phi_y$.
There is a partition of $y+V_1$ into $p$ affine subspaces of 
codimension 1, on each of which $\phi_y$ is constant, from
which it follows that there is a subspace $W_2$ of codimension 
$p^d$ on which \emph{all} of the functions $\phi_y$ are
constant.

Putting these two facts together, we obtain a subspace $V_2$ of codimension at most $r+p^d$ such that $Q$ is constant on cosets of $V_2$. 

Now let $V=V_0 \cap V_2$, so that $V$ has codimension at most $\d^{-4}T^4+r+p^d$.
Since $V$ is a subspace of $V_0$, we have that $\|U-U*\mu_V\|_2\leq\d$. Now $Q$ has constant modulus 1, so that $\|QU-Q(U*\mu_V)\|_2\leq\d$, and
since $Q$ is constant on cosets of $V$, $Q(U*\mu_V)=(QU)*\mu_V$. The result is proved.
\end{proof}

Finally, let us prove a couple more easy technical lemmas on how various $L_p$ and uniformity norms interact with the convolution operator. 


\begin{lemma} \label{shrink} 
Let $V$ be a subspace of $\F_p^n$ and let $f$ be a function from $\F_p^n$ to $\C$. Let $\|.\|$
be any translation-invariant norm defined on such functions. Then $\|f*\mu_V\|\leq\|f\|$.
\end{lemma}

\begin{proof}
If we write $f_v(x)$ for $f(x+v)$, we find that
$f*\mu_V=\E_{v \in V} f_v$ and we know that all the functions $f_v$ have
the same norm as $f$. The lemma therefore follows from the
triangle inequality.
\end{proof}


\begin{lemma}\label{U2vsL2}
Let $V$ be a linear subspace of codimension $r$ on $\F_p^n$ and 
let $f$ be a function from $\F_p^n$ to $\C$ that is constant on 
the cosets of $V$. Then $\|f\|_{U^2}\geq p^{-r/4}\|f\|_2$
and $\|f\|_{U^2}^*\leq p^{r/4}\|f\|_2$.
\end{lemma}

\begin{proof}
Let $\Gamma$ be a linear map from 
$\F_p^n$ to $\F_p^r$ with kernel $V$. Let $g:\F_p^r\rightarrow\C$
be defined by the formula $f(x)=g(\Gamma x)$, which is well-defined
since $f$ is constant on translates of $V$. It is easy to see that
$\|f\|_2=\|g\|_2$ and $\|f\|_{U^2}=\|g\|_{U^2}$. But
\begin{equation*}
\|g\|_{U^2}^4=\E_y|\E_x g(x)\overline{g(x+y)}|^2
\geq p^{-r}(\E_x|g(x)|^2)^2=p^{-r}\|g\|_2^4.
\end{equation*}
This proves the first part. For the second part, we know that $\|f\|_{U^2}^*$ is the maximum of $\langle f,h\rangle$ over all functions $h$ such that 
$\|h\|_{U^2}\leq 1$. Now replacing $h$ by $h*\mu_V$ does
not affect the inner product $\langle f,h\rangle$ and does not
increase $\|h\|_{U^2}$. Therefore, the maximum must be achieved
by a function $h$ that is constant on the cosets of $V$. 
But then $|\langle f,h\rangle|\leq\|f\|_2\|h\|_2$, which is
at most $p^{r/4}\|f\|_2\|h\|_{U^2}$, by the first part. This
completes the proof of the lemma.
\end{proof}

Before stating the main result of this section, let us recall Lemma 3.4 from our previous paper \cite{Gowers:2007tcs}, which says that the rank of a bilinear form cannot decrease too much when restricted to a smaller subspace. In that paper we gave a proof based on algebraic arguments. Here we shall present a different approach, included for the sake of completeness, which turns out to be more generalizable to higher-degree forms \cite{Gowers:2009lfuII} as well as locally-defined quadratic forms \cite{Gowers:2009lfuIII}.


\begin{lemma}\label{rankrestr}
Let $\beta$ be a symmetric bilinear form of rank $r$ on $\F_p^n$ 
and let $W$ be a subspace of $\F_p^n$ of codimension $d$. Then 
the rank of the restriction of $\beta$ to $W$ is at least $r-2d$.
\end{lemma}

\begin{proof} 
We first observe that for any fixed $x \in \F_{p}^{n}$, $\E_{y \in \F_{p}^{n}} \omega^{\beta(x,y)}$ is either 0 or 1, and hence 
\[p^{-r}=\E_{x,y \in \F_{p}^{n}} \omega^{\beta(x,y)}\geq \frac{|W|}{p^{n}}\E_{x \in W,y \in \F_{p}^{n}} \omega^{\beta(x,y)},\]
and similarly, for any fixed $y \in \F_{p}^{n}$, $\E_{x \in W} \omega^{\beta(x,y)}$ is either 0 or 1. It follows that 
\[\frac{|W|}{p^{n}}\E_{x \in W,y \in \F_{p}^{n}} \omega^{\beta(x,y)} \geq \frac{|W|^{2}}{p^{2n}}\E_{x,y \in W} \omega^{\beta(x,y)}=p^{-2d} p^{-r_{W}} ,\]
where we have written $r_{W}$ for the rank of the restriction of $\beta$ to $W$, and hence $r_{W} \geq r -2d$.
\end{proof}

Putting these results together, we obtain the main decomposition theorem that we
shall apply later to count certain types of linear configurations in uniform sets. It is similar to Theorem \ref{quadfouriersmallk},
but with the additional hypothesis that $f$ is highly uniform, which allows us to draw the stronger conclusion that all the quadratic averages used have
high rank.


\begin{theorem}\label{unifquadfourier}
Let $f:\F_p^n\rightarrow\C$ be a function such that 
$\|f\|_2\leq 1$. Then for every $\delta>0$ there exists
a constant $C$ such that for every $R_0$ there exists a
constant $c$ with the following property.

Let $d=(2/\d)^{C_p}$ and $C=(2/\d^2)^{C_p}$. Suppose that $\|f\|_{U^2}\leq c$. Then $f$ has a decomposition of the form
\begin{equation*}
f(x)=\sm i k Q_i(x)U_i(x)+g(x)+h(x),
\end{equation*}
where $k\leq C^2/\d^2$, the $Q_i$ are quadratic averages on $\F_p^n$ of complexity at most $d$, $\sm i k\|U_i\|_{U^2}^*\leq 2^{3/2}C^{4}\d^{-3}p^{3d/4}$, 
$\sm i k\|U_i\|_\infty\leq C$, $\|g\|_1\leq 7\d$ and
$\|h\|_{U^3}\leq2\d$. In addition, each quadratic average $Q_i$ has rank at least $R_0$ provided that $c$ satisfies the inequality $c\leq p^{-p^{16d}R_0}$.
\end{theorem}

\begin{proof}
Let us begin by applying Theorem \ref{quadfouriersmallk} to obtain a decomposition of the form
\begin{equation*}
f(x)=\sm i k Q_i(x)U_i(x)+g'(x)+h'(x),
\end{equation*}
where $k\leq C^2/\d^2$, the $Q_i$ are quadratic averages on $\F_p^n$ of complexity at most $d$, $\sm i k\|U_i\|_{U^2}^*\leq
2^{3/2}C^{4}\d^{-3}p^{3d/4}$, $\sm i k\|U_i\|_\infty\leq C$,
$\|g'\|_1\leq 2\d$ and $\|h'\|_{U^3}\leq\d$.

Let us assume that the quadratic averages $Q_i$ are arranged in
increasing order of rank. Now let $m=k$, let $r$ be such that
$p^{-r/2}=C^{2}/\d^2$ and let $t=p^{8d}+2d+2r$. Applying Lemma
\ref{rankgap} we obtain positive integers $R\in[R_0,m^k(R_0+t)]$ and
$s\in\{0,1,2,\dots,k\}$ such that $Q_i$ has rank at most $R$ when
$i\leq s$ and has rank at least $mR+t$ when $i>s$. Let $f_L=\sm i
sQ_iU_i$ and $f_H=\sum_{i=s+1}^kQ_iU_i$.

Let $T=2^{3/2}C^{4}\d^{-3}p^{3d/4}$. Theorem \ref{quadfouriersmallk}
tells us that $\sum_i\|U_i\|_{U^2}^*$, and hence each individual
$\|U_i\|_{U^2}^*$, is at most $T$.

Let $\eta=\d/k$. By Corollary \ref{quadavetimesu2*}, for every $i\leq s$ there is a linear subspace $V_i$ of codimension at
most $\eta^{-4}T^4+R+p^d$ such that
$\|Q_iU_i-(Q_iU_i) *\mu_{V_i}\|_2\leq\eta$. Let $V$ be the intersection of all the subspaces $V_i$. Then $V$ has codimension at
most $k(\eta^{-4}T^4+R+p^d)\leq kR+p^{8d}$ for sufficiently small $\d$, and $\|f_L-f_L*\mu_V\|_2\leq
k\eta\leq\d$.

Now let us return to our decomposition $f=f_L+f_H+g'+h'$. We shall
convolve both sides with the measure $\mu_V$ of the subspace $V$ just constructed and consider its effect on the $L_2$ norms of $f$ and $f_H$. If both of these are small, it will allow us to approximate $f_L$ by a quadratically uniform function up to an error in $L_1$.

First, since $\|f\|_{U^2}\leq c$, Lemma \ref{shrink} implies that
$\|f*\mu_V\|_{U^2}\leq c$. It follows from Lemma \ref{U2vsL2} that
$\|f*\mu_V\|_2\leq cp^{(kR+p^{8d})/4}$, which is at most $\d$ by our
choice of $c$.
Next, let us look at $f_H$. For each $i>s$, the quadratic average
$Q_i$ has rank at least $kR+t$. Recall from the proof of Theorem \ref{quadfouriersmallk} that $Q_iU_i$ is equal to a sum
of the form $\sum_{j\in A_i}\lambda_jQ_j$ with $\sum_i\sum_{j\in
A_i}|\lambda_j|\leq C$. For each $j$, write $(V_j,q_j)$ for
the base of $Q_j$. Then we had the additional property that
for every $j\in A_i$ the rank of $q_i-q_j$, considered as a 
quadratic form on $V_i\cap V_j$, was at most the $r$ chosen 
earlier---that is, the $r$ such that
$p^{r/2}=2C/\d^2$.

If we consider $Q_i$ as a quadratic average with base $(V_i\cap V_j,q_i)$,
then by Lemma \ref{rankrestr} it has rank at least $kR+t-2d$, and $Q_j$, considered as a quadratic average with base 
$(V_i\cap V_j,q_j)$ therefore has rank at least $kR+t-2d-r$. It 
follows from Lemma \ref{quadaveu2} that 
$\|Q_j\|_{U^2}\leq p^{-(kR+t-2d-r)/4}$, and hence from the triangle
inequality that 
$\|f_H\|_{U^2}\leq\sum_{i>s}\|Q_iU_i\|_{U^2}\leq Cp^{-(kR+t-2d-r)/4}$. 
Therefore, by the same argument as we used for $f$, we find that
\[\|f_H*\mu_V\|_2\leq Cp^{-(kR+t-2d-r)/4}p^{(kR+p^{8d})/4}
=Cp^{(p^{8d}+2d+r-t)/4}=Cp^{-r/4}\]
By our choice of $r$, this is less than $\d$. 

As for $g'$ and $h'$, we know by Lemma \ref{shrink} that
$\|g'*\mu_V\|_1\leq\|g'\|_1\leq 2\d$ and $\|h'*\mu_V\|_{U^3}\leq\d$.
Since $\|f_L-f_L*\mu_V\|_2\leq\d$, we have ended up showing that
$f_L$ can be written as a sum $g''+h''$, where $\|g''\|_1\leq 5\d$
and $\|h''\|_{U^3}\leq\d$. Therefore, we can write $f=f_H+g+h$
with $\|g\|_1\leq 7\d$, $\|h\|_{U^3}\leq 2\d$. This proves
the theorem.
\end{proof}

\section{Proof of Theorem \ref{impbound}}

The aim of this section is to show that if $f$ is a function of
the form $\sum_j U_j Q_j$, where the $Q_j$ are high-rank
quadratic averages and $\sum_j \|U_j\|_\infty$ is not too large,
then $\E_{\x\in(\F_p^n)^d}\prod_i f(L_i(\x))$
is small in modulus whenever the linear forms $L_1,\dots,L_m$ are square independent. It is then relatively straightforward to deduce Theorem \ref{impbound}.

First, let us prove some lemmas that will help us use
the high-rank condition on the quadratic averages in conjunction with the  
square independence of the linear forms. The first result states that if the bilinear form $\beta$ has high rank, then the phase function $\omega^{\beta(x,y)}$
is quasirandom.


\begin{lemma}\label{bilinearqr}
Let $\beta$ be a bilinear form of rank at least $r$ on a subspace $V$ of $\F_p^{n}$, and let
$g$ and $h$ be two functions with $\|g\|_\infty$ and
$\|h\|_\infty$ at most 1. Then 
\[|\E_{x,y}\omega^{\beta(x,y)}g(x)h(y)|\leq p^{-r/2}.\]
\end{lemma}

\begin{proof}
This lemma can be proved either directly (as we shall do) or indirectly, by first estimating the rectangle norm of the function and applying standard results in the theory of quasirandomness. Either way, the proof is a standard application of the Cauchy-Schwarz inequality.
\[|\E_{x,y}\omega^{\beta(x,y)}g(x)h(y)|^2\leq
\E_x|\E_y\omega^{\beta(x,y)}g(x)h(y)|^2
\leq \E_x|\E_y\omega^{\beta(x,y)}h(y)|^2.\]

The latter expression can be expanded as
\[\E_{y,y'}h(y)h(y')\E_x\omega^{\beta(x,y-y')}
\leq\E_{y,y'}|\E_x\omega^{\beta(x,y-y')}|.\]

Now $\beta(x,y-y')$ depends linearly on $x$, so 
$\E_x\omega^{\beta(x,y-y')}$ is zero unless 
$\omega^{\beta(x,y-y')}$ is constant. That is, $\E_x\omega^{\beta(x,y-y')}$ is zero unless 
$y-y'$ belongs to the kernel of $\beta'$. Otherwise, it has modulus 1. Since $\beta$ has 
rank at least $r$, the probability, for each $y$, that $y-y'$ 
belongs to the annihilator is at most $p^{-r}$. Therefore, 
$$\E_{y,y'}|\E_x\omega^{\beta(x,y-y')}|\leq p^{-r}.$$
The result follows on taking square roots.
\end{proof}

In order to show that the average over a product of quadratic phases is small, we actually only need one of the quadratic phases involved to have high rank.


\begin{lemma} \label{justonebilinear}
Let $d$ be a positive integer and for every pair 
$(u,v)\in[d]^2$ let $\beta_{uv}$ be a bilinear
form on $\F_p^n$ taking variables $x_{u}$ and $x_{v}$. For $u \in [d]$, let $\phi_{u}$ be a linear functional on $(\F_p^n)^d$ in the variable $x_{u}$. Suppose 
that the rank of $\beta_{uv}$ is at least $r$ for at least one pair $(u,v)$. Then
\[
\Bigl|\E_{\x\in (\F_p^n)^d}\omega^{\sum_{u,v}\beta_{uv}(x_u,x_v)+\sum_{u}\phi_{u}(\x_{u})}\Bigr|
\leq p^{-r/2}.
\]
\end{lemma}

\begin{proof}
Let us assume first that $\beta_{uu}$ has rank at least $r$ for
some $u$. If we fix the values of $x_v$ for every $v\ne u$, then
the sum in the exponent takes the form $\beta_{uu}(x_u,x_u)+\gamma(x_u)$
for some linear functional $\gamma$. Therefore, by Lemma \ref{gauss} the
expectation over $x_u$ has modulus at most $p^{-r/2}$. Since this 
is true for every choice of the other $x_v$, the whole expectation 
has modulus at most $p^{-r/2}$. 

Now let us assume that $\beta_{uv}$ has rank at least $r$ for
some pair $(u,v)$ with $u\ne v$. This time, let us fix all the
variables apart from $x_u$ and $x_v$. Now the sum in the exponent
takes the form 
$$2\beta_{uv}(x_u,x_v)+\phi(x_u)+\psi(x_v)$$
so by Lemma \ref{bilinearqr} the expectation over $x_u$ and
$x_v$ is at most $p^{-r/2}$. Again, since this is true for every
possible choice of the other variables, the whole expectation
is at most $p^{-r/2}$.
\end{proof}

Next we shall show that if we have a set of bilinear forms of high rank, then at least one of the linear combinations arising from a square-independent system $L_1, L_2, \dots, L_m$ must have fairly high rank.


\begin{lemma}\label{rankaverage}
Let $V$ be a subspace of $\F_p^{n}$ and let $\beta_1, \dots,
\beta_m$ be bilinear forms on $V$ with rank at least $r$. Let $B$ be
an invertible $m\times m$ matrix with entries $b_{ij} \in \F_p$. Then
at least one of the bilinear forms $\eta_j=\sum_{i=1}^m b_{ij}\beta_i$
has rank at least $r/m$.
\end{lemma}

\begin{proof}
It follows from the assumption that $B$ is invertible that
$\beta=B^{-1}\eta$. But the rank of a linear combination of the
$\eta_i$ is at most the sum of the ranks of the $\eta_i$.
\end{proof}


\begin{corollary}\label{rankaveragecor}
Suppose that $L_i(\x)=\sum_{u=1}^d c_{iu} x_u$, $i=1, 2, \dots, m$, is a
square-independent system. Suppose that each of the (not necessarily
distinct) bilinear forms $\beta_i$, $i=1, 2, \dots, m$, has rank at least
$r$. Then at least one of the bilinear forms $\beta_{uv}=\sum_{i=1}^m
c_{iu}c_{iv} \beta_i$ has rank at least $r/m$.
\end{corollary}

\begin{proof}
For each $i=1, 2, \dots, m$, let $M_i$ be the matrix $(c_{iu}
c_{iv})_{u,v}$. Square independence implies that the matrices $M_i$
are linearly independent over $\F_p$. This implies that the rank of
the $d^2 \times m$ matrix whose $((u,v),i)$ entry is $c_{iu}c_{iv}$ is
$m$. The rows of this matrix are the $(d\times d)$ matrices $M_1, \dots, M_m$. The
columns are the vectors $C_{uv}=(c_{1u}c_{1v}, c_{2u}c_{2v}, \dots,
c_{mu}c_{mv})$. Since row rank equals column rank, we can find $m$
linearly independent vectors $C_{uv}$. Now apply Lemma
\ref{rankaverage} to the bilinear forms $\beta_{uv} = \sum_{i=1}^m
(C_{uv})_{i} \beta_i$ to obtain the result.
\end{proof}

We are now in a position to prove the key ingredient of Theorem \ref{impbound}.


\begin{proposition}\label{quadraticpart}
Let $C$, $D$, $R$, $T$ and $\delta$ be positive constants. For each $i=1,2,\dots, m$, let $f_{i}=\sum_{j=1}^{k_{i}} U_j^{(i)}Q_j^{(i)}$ be a linear combination of quadratic averages $Q_j^{(i)}$ on $\F_p^n$, each of
rank at least $R$ and complexity at most $D$, such that
$\sum_{j=1}^{k_{i}} \|U_j^{(i)}\|_\infty\leq C$ and $\sum_{j=1}^{k_{i}}\|U_j^{(i)}\|_{U^2}^*\leq T$. Let $d$ and $m$ be positive integers, and
let $\seq L m$ be a square-independent system of $m$ linear forms in $d$
variables. Then 
\[\left|\E_{\x\in(\F_p^n)^d}\prod_if_i(L_i(\x))\right|\leq (C^m p^{D+\delta^{-4}T^4-R/2m}+ \delta)\prod_{i=1}^{m}k_{i}.\]
\end{proposition}

\begin{proof}
Since $f$ is a sum of quadratic averages, the expectation in question 
can be split up into a sum of terms of the form
\[\E_{\x\in(\F_p^n)^d}\prod_{i=1}^{m}(U_{j_i}^{(i)}Q_{j_i}^{(i)})(L_i(\x)).\]
Let us obtain an upper bound for the size of one of these terms.
For ease of notation, we shall take the sequence $(\seq j m) \in [k_{1}]\times[k_{2}]\times\dots\times[k_{m}]$
to be the sequence $(1,2,\dots,m)$. That is, we let $\seq Q m$
be an arbitrary sequence of quadratic averages of rank at least
$R$ and complexity at most $D$, and we let $U_i$ an arbitrary sequence of bounded anti-uniform functions. We shall obtain an upper bound for
the modulus of $\E_{\x\in(\F_p^n)^d}\prod_i (U_iQ_i)(L_i(\x))$.

By Lemma \ref{bogolyubov}, we can find, for each $i=1, 2, \dots, m$, a subspace $W_i$ of codimension at most $\delta^{-4}T^4$ such that $\|U_i -U_i*\mu_{W_i}\|_2 \leq \delta$, in other words, each $U_i$ can be approximated by a function that is constant on translates of $W_i$ and still bounded. Set $W=W_1 \cap \dots \cap W_m$, which is a subspace of codimension at most $m\delta^{-4}T^4$. Then 
$\|U_i -U_i*\mu_{W}\|_2 \leq \delta$ for all $i=1, 2, \dots, m$. We shall replace $U_i$ by $U_{i}*\mu_{W}$ in the above average, incurring an error of at most $\delta$, and from now on focus on $\E_{\x\in(\F_p^n)^d}\prod_i ((U_i*\mu_W)Q_i)(L_i(\x))$.

For each $i$, let us write $(V_i,q_i)$ for the base of the quadratic 
average $Q_i$, and let $V=W\cap V_1\cap\dots\cap V_m$. Then we can, if we 
wish, regard each $Q_i$ as a quadratic average with base $(V,q_i)$. 
Since each $V_i$ has codimension at most $D$ in $\F_p^n$ and $W$ had codimension at most $m\delta^{-4}T^4$, the
codimension of $V$ is at most $m(D+\delta^{-4}T^4)$, from which it follows that 
the rank of $Q_i$, when considered in this new way, is at least
$R-2m(D+\delta^{-4}T^4)$ by Lemma \ref{rankrestr}.

We now split the expectation $\E_{\x\in(\F_p^n)^d}\prod_i ((U_i*\mu_W)Q_i)(L_i(\x))$ even further, according to the particular set of translates of $V$
that the components $\seq x d$ of $\x$ belong to. Let $\seq V d$
be arbitrary translates of $V$ and let us obtain a bound for the size of
\[\E_{x_1\in V_1}\E_{x_2\in V_2}\dots\E_{x_d\in V_d}\prod_{i=1}^{m}((U_i*\mu_W)Q_i)(L_i(\x)).\]
Now if each $x_j$ is confined to $V_j$, then $L_i(\x)$ is confined to some
particular translate $y+V$ of $V$. On this translate, $U_i*\mu_W$ is constant, say equal to $\lambda_y$, and  $Q_i(x)$ is given by a formula of the form $\omega^{q_i(x-y)+\phi_y(x-y)}$,
where $\phi$ is linear. It follows that the expectation we are trying
to estimate is equal to a quantity of the form
$\E_{\x\in V^d}\prod_i\omega^{q_i(L_i(\x))+\phi_i(L_i(\x))}$,
where each $q_i$ has rank at least $R-2m(D+\delta^{-4}T^4)$, and we temporarily disregard the product of the coefficients $\prod_{i=1}^m\lambda_{y_i}$.

Now for each $i=1,2,\dots, m$, let the linear form $L_i(\x)$ be given by the formula $\sum_{u=1}^d c_{iu}x_u$. Then, writing $\b_i$ for the bilinear form associated with $q_i$, we have 
\[\sm i mq_i(L_i(\x))=\sum_{u,v=1}^{d}\sm i mc_{iu}c_{iv}\b_i(x_u,x_v).\]
For each $u$ and $v$, let us write $\b_{uv}$ for the bilinear form
$\sm i mc_{iu}c_{iv}\b_i$. Lemma \ref{rankaveragecor} implies
that at least one of the bilinear forms $\b_{uv}$ has rank at
least $R/m-2(D+\delta^{-4}T^4)$.

Write $\sm i m\phi_i(L_{i}(\x))=\sum_{i=1}^{m}\phi_{i}(\sum_{u=1}^{d}c_{iu}x_{u})=\sum_{u=1}^{d}\sum_{i=1}^{m}c_{iu}\phi_{i}(x_{u})=\sum_{u=1}^{d}\phi_{u}(x_{u})$. Then
\[\E_{\x\in V^d}\prod_{i=1}^{m}\omega^{q_i(L_i(\x))+\phi_i(L_i(\x))}=\E_{\x\in V^d}\omega^{\sum_{u,v}\b_{uv}(x_u,x_v)+\sum_{u}\phi_{u}(\x_{u})}.\]
Applying Lemma \ref{justonebilinear}, we may deduce that
\[|\E_{\x\in V^d}\omega^{\sum_{u,v}\b_{uv}(x_u,x_v)+\phi(\x)}| \leq p^{(D+\delta^{-4}T^4)-R/2m}\]
Since this estimate did not depend on our choice of translates of $V$, it follows that 
\[|\E_{\x\in(\F_p^n)^d}\prod_i ((U_i*\mu_W)Q_i)(L_i(\x))|\leq C^m p^{(D+\delta^{-4}T^4)-R/2m},\] 
where the sum over the product of the constant coefficients $\lambda_i$, which we had temporarily neglected, contributes the additional factor of $C^m$.
Since this was true for an arbitrary choice of quadratic averages of rank at least $R$, the statement of the lemma follows.
\end{proof}

Before we are able to prove Theorem \ref{impbound} we need one more standard result that will allow us to neglect the quadratically uniform part of the decomposition. The following statement is implicit in Green and Tao \cite{Green:2006lep}, and was also a major ingredient in \cite{Gowers:2007tcs}. The proof is a repeated application of the Cauchy-Schwarz inequality together with a suitable reparametrization of the linear system under consideration.


\begin{theorem}\label{gvnmod}
Let $f_1,\dots,f_m$ be functions $\F_p^n$, and let $\L_{1}, L_{2}, \dots, L_{m}$ be a
linear system of Cauchy-Schwarz complexity $k$ consisting of $m$ forms in $d$
variables. Then
\[\left|\E_{\x \in (\F_p^n)^d} \prod_{i=1}^m f_i(L_i(\x)) \right|\leq
\min_i\|f_i\|_{U^{k+1}}\prod_{j\neq i} \|f_j\|_\infty.\]
\end{theorem}

Let us now put all the technical results from the preceding two sections together to give an improved bound for Theorem \ref{oldresult}, and thus prove Theorem \ref{impbound}.


\begin{proof}[Proof of Theorem \ref{impbound}]
Let $\eps>0$, and let $c>0$ be chosen in terms of $\eps$ later. Given $f:\F_p^n\rightarrow[-1,1]$ with $\|f\|_{U^2} \leq
c$ we first apply Theorem \ref{unifquadfourier} with
$\delta_1=\eps/(18m)$ to obtain a decomposition
$$f=f_1+g_1+h_1,$$ where $f_1=\sum_jU_j^{(1)} Q_j^{(1)}$ with
$\sum_j\|U_j^{(1)}\|_\infty\leq M_1$, $\|g_1\|_1 \leq 7 \delta_1$ and
$\|h_1\|_{U^3} \leq 2 \delta_1$. We have carefully ensured that each
quadratic average $Q_j^{(1)}$ has rank at last $R_0$ for some $R_0$ to be chosen
later, and $M_1$ is a function of $\delta_1$ only, which can be taken to
equal $(2\delta_1^{-2})^{C_p}$.  Recall that we want to show that
$$\E_{\x\in(\F_p^n)^d}\prod_{i=1}^m f(L_i(\x))$$ 
is bounded in absolute value by $\eps$ for sufficiently uniform $f$. We begin by replacing
the first $f$ in the product by $g_1+h_1$. The product involving $g_1$
yields an error term of $7 \delta_1$ since all the remaining factors
have $L_{\infty}$ norm bounded by $1$, while the product involving
$h_1$ yields an error of $2\delta_1$ by Theorem \ref{gvnmod} above with $k=2$. Our choice of $\delta_1$ implies that the sum of these two errors is at most $\eps/(2m)$.

Now we apply Theorem \ref{unifquadfourier} again, this time with
$\delta_2=\eps/(18mM_1)$, to obtain a decomposition
$$f=f_2+g_2+h_2,$$ where $f_2=\sum_jU_j^{(2)} Q_j^{(2)}$ with
$\sum_j\|U_j^{(2)}\|_\infty\leq M_2$, $\|g_2\|_1 \leq 7 \delta_2$ and
$\|h_2\|_{U^3} \leq 2 \delta_2$. When replacing the first instance of $f$ in
$$\E_{\x\in(\F_p^n)^d}f_1(L_1(\x))\prod_{i=2}^{m} f(L_i(\x))$$ 
with $g_2+h_2$, the product involving $g_2$ now contributes an
error term of at most $7\delta_2 M_1$ (since $\|f_1\|_{\infty}\leq
M_1$). By Theorem \ref{gvnmod} it follows that the contribution from the product involving $h_2$ is
bounded above by $2\delta_2 M_1$. Therefore the total error incurred
is at most $9\delta_2 M_1$, which is at most $\eps/(2m)$ by our
choice of $\delta_2$.

When we come to apply Theorem \ref{unifquadfourier} to the $k$th instance of
$f$ in the product, we need to do so with $\delta_k$ satisfying $9
\delta_k M_1 \dots M_{k-1} \leq \eps/(2m)$ for $k=2,\dots,m$. This
ensures that up to an error of $\eps/2$, it suffices to consider the product
$$\E_{\x\in(\F_p^n)^d}\prod_{i=1}^m f_i(L_i(\x)).$$
Since each $M_k$ is a polynomial in $\delta_k^{-1}$, and since $\delta_1$ was chosen proportional to $\eps$, it is easy to see that $M_m$ will be bounded above by a polynomial of $\eps^{-1}$. In fact, it is not difficult to establish that the bound on $M_m$ will be of the form $c_{m,p}\eps^{-(4C_p)^m}$, where $c_{m,p}$ will be a constant depending on $m$ and $p$ only.

Recall that each $f_i$ was of the form $\sum_{j=1}^{k_i} U_j^{(i)} Q_j^{(i)}$ with
$\sum_j\|U_j^{(i)} \|_\infty\leq M_i$, $\sum_j\|U_j^{(i)} \|_{U^2}^*\leq T_i$  and each $Q_j^{(i)}$ had rank at least $R_0$ and complexity at most $d_i$.
It is clear from the procedure we have applied that the parameters $d_i$, $k_i$, $M_i$ and $T_i$ are strictly increasing in $i$, and that we can take $d_m=(2/\delta_m)^{C_p}$, $k_m=(2/\delta_m^2)^{C_p}/\delta_m^2$, $M_m=(2/\delta_m^2)^{C_p}$ and $T_m=2^{3/2}M_m^4\delta_m^{-3}p^{3d_m/4}$ by Theorem \ref{unifquadfourier}.

Key Proposition \ref{quadraticpart} with $\delta= \eps k_m^{-m}/4$ now implies that 
$$\E_{\x\in(\F_p^n)^d}\prod_{i=1}^m f_i(L_i(\x))$$
is bounded in modulus by 
\[k_m^m M_m^m p^{d_m+\delta_{m}^{-4}(2^6M_m^{16}\delta_m^{-12}p^{3d_m})-R_0/2m}+\eps/4.\] 
Let us analyse this expression. As we have already remarked, at each step $\delta_l$ is a polynomial in $\eps$, and thus the quantities $d_m$, $k_m$ and $M_m$ are polynomial in $\eps^{-1}$. It follows that $R_0$ needs to be taken exponential in $\eps^{-1}$ at each step.

More precisely, $\delta_m$ can be chosen of the form $c_{m,p} \eps^{(4C_p)^{m-1}}$ for a constant $c_{m,p}$ only depending on $m$ and $p$, and hence $d_m$ can be assumed to be at most $c'_{m,p} \eps^{-(4C_p)^{m}}$. It therefore suffices to take $R_0$ of the form $\exp(c''_{m,p}\eps^{-(4C_p)^{m}})$.

However, in order to be able to choose $R_0$ as the minimum rank in each decomposition of $f$, we needed the uniformity parameter $c$ to satisfy $p^{-p^{16d_m}R_0}$, which is a function of the form $\exp(-\exp(c'''_{m,p}\eps^{-(4C_p)^{m}}))$.
\end{proof}

\section{Remarks}

An obvious question to ask is whether the bounds in Theorem \ref{impbound} can be improved further. We remark that it is possible to obtain a single exponential in Theorem \ref{impbound} if one works under the assumption of the so-called \emph{Polynomial Freiman-Ruzsa Conjecture}. This conjecture asserts that a subset $A$ of doubling $K$ can be covered by at most $C_1(K)$ translates of a subspace of size at most $C_2(K)|A|$, where both $C_1(K)$ and $C_2(K)$ are polynomial in $K$ (see for example \cite{Green:2005ffm}). It has recently been shown to be equivalent to polynomial bounds in Theorem \ref{u3inverse} \cite{Green:2009ebi}. Applying a local version of this conjecture, we find that $\delta_l$ and $M_l$ in the proof above still grow polynomially in $\eps^{-1}$, while the dimension $d$ remains logarithmic in $\eps^{-1}$, reducing the final bound to a single exponential. In the other direction, we do not know of a lower bound that is better than a power.


\end{document}